\documentclass{gtpart}     
\usepackage{pinlabel}  
\usepackage{graphicx}  
\usepackage[all]{xy}
\usepackage{amscd}
\usepackage{tikz}
\title{On the Tate spectrum of $\mathrm{tmf}$ at the prime 2}

\author{Scott M. Bailey}
\givenname{Scott}
\surname{Bailey}
\address{Clayton State University, Department of Mathematics, 2000 Clayton State Blvd., Morrow, GA 30260}
\email{ScottBailey@clayton.edu}
\urladdr{http://faculty.clayton.edu/sbailey21}

\author{Nicolas Ricka}
\givenname{Nicolas}
\surname{Ricka}
\address{Wayne State University, Department of Mathematics, 656 W. Kirby, Detroit, MI 48202}
\email{nicolas.ricka@wayne.edu}
\urladdr{https://sites.google.com/site/nicoricka/}
\keyword{Tate spectrum}
\subject{primary}{msc2010}{55P42}
\subject{secondary}{msc2010}{55T15}

\arxivreference{}
\arxivpassword{}

\newtheorem{thm}{Theorem}[section]    
\newtheorem{lem}[thm]{Lemma}          
\newtheorem{prop}[thm]{Proposition} 
\theoremstyle{definition}
\newtheorem{defn}[thm]{Definition}    
\newcommand{\tmf}{\mathrm{tmf}}
\newcommand{\ko}{\mathrm{kO}}
\newcommand{\SA}{\mathcal{A}}

\newcommand{\RR}{\mathbb{R}}
\newcommand{\ZZ}{\mathbb{Z}}
\newcommand{\FF}{\mathbb{F}}

\newcommand{\sm}{\wedge}
\newcommand{\sur}{/ \hspace{-0.07cm}/}
\newcommand{\Ext}{\operatorname{Ext}}

\newcommand{\sus}[1]{\Sigma^{#1}}
\newcommand{\dsus}[1]{\Omega^{#1}}

\begin{document}

\begin{abstract}    
Computations involving the root invariant prompted Mahowald and Shick to develop the slogan: ``the root invariant of $v_n$-periodic homotopy is $v_n$-torsion.'' While neither a proof, nor a precise statement, of this slogan appears in the literature, numerous authors have offered computational evidence in support of its fundamental idea. The root invariant is closely related to Mahowald's inverse limit description of the Tate spectrum, and computations have shown the Tate spectrum of $v_n$-periodic cohomology theories to be $v_n$-torsion. The purpose of this paper is to split the Tate spectrum of $\tmf$ as a wedge of suspensions of $\ko$, providing yet another example in support of the slogan to the existing literature.
\end{abstract}

\maketitle


\section{Introduction} \label{sec:introduction}

Let $\lambda$ denote the canonical line bundle over $\RR P^\infty = B(\ZZ/2\ZZ)$. For $\ell \in \ZZ$, define $P_\ell$ to be the Thom spectrum of $\ell\lambda$. Induced maps on the level of Thom spectra give a naturally defined inverse system of projective spaces

\begin{equation}\label{eq:inversesystem_projectivespaces}
\dots \to P_{n-1} \xrightarrow{\jmath_{n-1}} P_n \xrightarrow{\jmath_n} P_{n+1} \to \dots
\end{equation}

W.H. Lin \cite{Lin80} demonstrated that the homotopy limit of \eqref{eq:inversesystem_projectivespaces} has the homotopy type of a desuspended 2-complete sphere, i.e., 

\begin{equation} \label{eq:linstheorem}
\varprojlim P_n \simeq \hat{S}^{-1}
\end{equation}

Suppose $X$ is a finite complex and consider a cohomotopy class $\alpha \in [X,S^0]_j$. The equivalence \eqref{eq:linstheorem} guarantees the existence of a largest $\ell \in \ZZ$ such that the composite $\sus{j-1}X \xrightarrow{\alpha} S^{-1} \to P_\ell$ is nontrivial.  In particular, this induces a map $R(\alpha) : \sus{j-1} X \to S^\ell$, the homotopy class of which is called the root invariant of $\alpha$.  Computations inside the EHP sequence led Mahowald and Ravenel \cite{MR87} to conjecture that the root invariant of a $v_n$-periodic element is $v_{n+1}$-periodic. This prompted Mahowald and Shick \cite{MS88} to discuss the related slogan: ``The root invariant of $v_n$-periodic homotopy is $v_n$-torsion.'' They show, for a finite complex $X$ having a $v_n$-self map, that

\begin{equation} \label{eq:vntorsion}
\varprojlim \left([X,P_\ell](v_n^{-1})\right) = 0
\end{equation}

In particular, Mahowald and Shick point out that if $\alpha \in [X, S^0]$ is $v_n$-periodic then its root invariant, at least when considered as an element of $[X, P_\ell]$, is $v_n$-torsion.

Let $X$ be a spectrum. Mahowald's description of the Tate spectrum of $X$ is the homotopy inverse limit
\begin{equation}\label{eq:tatespectrum_invlimconstruction}
t(X) = \Sigma \varprojlim (P_\ell \sm X)
\end{equation}
If $X$ is a finite spectrum, then \eqref{eq:tatespectrum_invlimconstruction} implies the Tate spectrum functor corresponds to completion at the prime 2. This is certainly not the case for all spectra since homotopy limits do not commute with the smash product. 

While neither a proof, nor even a precise statement of the phenomenon suggested by the slogan has appeared in the literature, many authors have demonstrated the Tate spectrum functor sends $v_n$-periodic cohomology theories to $v_{n-1}$-periodic theories:

\begin{description}
	\item[1984:] Davis and Mahowald \cite{DM84}, for $p=2$, show $t(\ko) \simeq \bigvee_{j \in \ZZ} \sus{4j} \widehat{H\ZZ}$;
	\item[1986:] Davis, Johnson, Klippenstein, Mahowald, and Wegmann \cite{DM86} demonstrate that, if $p$ any prime and $q = 2(p-1)$, then there are equivalences of $p$-complete spectra $t(BP\langle 2 \rangle) \simeq \prod_{j\in\ZZ} \sus{qj} \widehat{BP}\langle 1 \rangle$, and conjecture a similar splitting of $t(BP\langle n \rangle)$;
	\item[1998:] Ando, Morava, and Sadofsky \cite{AMS98} prove the existence of a ring isomorphism $(tE(n)_*)^{\wedge}_{I_{n-1}} \cong E(n-1)_*((x))^\wedge_{I_{n-1}}$ where $I_{n-1} = (p,v_1,\hdots,v_{n-2})$ and construct a map of spectra $\bigvee_{j \in \ZZ} \sus{2j} E(n-1) \to tE(n)^\wedge_{I_{n-1}}$ which, after completion at $I_{n-1}$ (or equivalently after localization with respect to the $(n-1)$st Morava $K$-theory) induces the isomorphism of homotopy groups.
\end{description}

The purpose of this paper is to provide yet another example to the literature. Let $\tmf$ denote the connective ring spectrum of topological modular forms (see \cite{Beh07, Goe09, Lur07}) at the prime 2. The main theorem is:

\begin{thm} \label{thm:maintheorem}
	There is a weak equivalence of spectra 
\begin{equation}\label{eq:maintheorem}
t(\tmf) \simeq  \prod_{i \in \Z} \Sigma^{8i} \ko,
\end{equation}
\end{thm}

In the context of the above machinery, computations involving the homotopy of $t(\tmf)$ greatly benefit from Mahowald's inverse limit description of the Tate spectrum.  However, the Tate spectrum functor conserves other properties, such as ring structure, of the spectrum.  This fact, however, is not immediately clear from the inverse limit point of view. On the other hand, such a structure is clear when placed in the framework established by Greenlees and May \cite{GM95a}. In their notation: let $G$ be a compact Lie group, $EG$ a free contractible $G$-space and $\tilde{E}G$ the cofiber of the map $EG_+ \to S^0$. If $k_G$ is a $G$-spectrum, then $\mathbf{t}(k_G) = F(EG_+,k_G) \sm \tilde{E}G$, where $F(EG_+,k_G)$ is the function $G$-spectrum of maps $EG_+ \to k_G$, is the Tate spectrum of $k_G$. Since $EG_+$ is equipped with a coproduct, if $k_G$ is a ring spectrum then $F(EG_+, k_G)$ is also a ring spectrum. Combining this with the product on $\tilde{E}G$, $\mathbf{t}(k_G)$ is also a ring. Lewis-May fixed points give a lax monoidal functor, so $\mathbf{t}(k_G)^G$ also has a ring structure.

The link to Mahowald's inverse limit description is as follows \cite{GM95a}: If $G$ is cyclic of order 2 and $k_G$ is the equivariant $G$-spectrum associated to a non-equivariant spectrum $k$, then there is a homotopy equivalence $\mathbf{t}(k_G)^G \simeq \Sigma \varprojlim (P_\ell \sm k) = t(k)$.

With the above correspondence in mind, we can restate Theorem~\ref{thm:maintheorem} as
\begin{thm}\label{thm:maintheorem_restated}
	There is a weak equivalence of ring spectra
	\begin{equation}\label{eq:maintheorem_restated}
		t(\tmf) \simeq \ko[x^{\pm 1}]
	\end{equation}
	where $x$ is degree 8.
\end{thm}

\section{Some particular $\SA(2)$-modules} \label{sec:A2modules}

Let $\SA$ denote the mod-2 Steenrod algebra generated by the squaring operations $\{ Sq^{2^i} \}_{i\geq 0}$. Let $M$ be an $\SA$-module, and consider the $\SA$-modules $\SA\sur\SA(n) \otimes M$ via the diagonal action and $\SA \otimes_{\SA(n)} M$ via left action on $\SA$. There is an isomorphism

\begin{equation}\label{eq:SAnisomorphism}
\Phi : \SA\sur\SA(n) \otimes M \to \SA \otimes_{\SA(n)} M
\end{equation}  

defined by $\Phi(a \otimes m) = \sum a' \otimes a'' m$ where $\psi(a) = \sum a' \otimes a''$ is the coproduct on $\SA$.  This isomorphism induces a change-of-rings isomorphism on the level of $\Ext$-groups

\begin{equation}\label{eq:changeofrings}
\Ext_\SA^{s,t}(\SA\sur\SA(n) \otimes M, N) \cong \Ext_\SA^{s,t}(\SA\otimes_{\SA(n)} M, N) \cong \Ext_{\SA(n)}^{s,t}(M,N)
\end{equation}

which is often invoked to simplify computations within the Adams sepctral sequence $E_2$-term.  For instance, since $H^*(\tmf) \cong \SA\sur\SA(2)$ \cite{Mat13, HM98}, to compute the homotopy groups of $P_i \sm \tmf$ it suffices to understand the left $\SA(2)$-module structure of $H^*P_i$. 

The following two propositions are results of Lin, Davis, Mahowald, and Adams \cite{LDM80}.
 
\begin{prop}\label{prop:SAPi}
	As a $\FF_2$-vector space, $H^* P_i = x^i \FF_2[x]$, where $x$ is in degree 1. The action of $\SA$ on $H^*P_i$ is determined by
	\begin{equation} \label{eq:SAPi}
		Sq^j x^k = \binom{k}{j} x^{j+k}
	\end{equation}
\end{prop}

Denote by $\FF_2 [x^{\pm 1}]$ the colimit of these $\SA$-modules. Note that this is certainly not the cohomology of the limit given by \eqref{eq:linstheorem}. A consequence is the following interpretation of the $\SA(2)$-module structure of some quotients of $\FF_2[x^{\pm1}]$:

\begin{prop} \label{prop:quotientSA2}
	Let $F_\ell$ be the sub-$\SA(2)$-module of $\FF_2[x^{\pm1}]$ generated by the classes in degree less than $\ell$. Then there is an isomorphism of $\SA(2)$-modules
	\begin{equation}\label{eq:quotientSA2}
		\FF_2[x^{\pm1}]/F_\ell  \cong \bigoplus_{j \geq \frac{\ell+1}{8}} \sus{8j-1}\SA(2)\sur\SA(1)
	\end{equation} 
\end{prop}

\begin{defn}
	Let $L_0$ denote the spectrum $(S^1 \cup_2 e^2 \cup_\eta e^4 \cup_\nu e^8)_+$. By construction, $H^*L_0 = \FF_2\{ 1, x, Sq^1x, Sq^2Sq^1x, Sq^4Sq^2Sq^1 x \}$, with the action of $\SA$ indicated by the names of the elements.
\end{defn}

\begin{prop} \label{prop:associategradedP-1}
	There is a filtration of $\SA(2)$-modules of $H^* P_{-1}$ with associated graded $H^*L_0 \oplus \bigoplus_{j \geq 0} \sus{8j-1} \SA(2)\sur\SA(1)$.
\end{prop}

\begin{proof}
Note that there is an isomorphism of $\SA(2)$-modules $\FF_2[x^{\pm1}] \cong \left(\FF_2[x^{\pm1}]/F_{-1}\right) \oplus F_{-1}$ so that $H^*P_{-1} \cong \left( \FF_2[x^{\pm1}]/F_{-1}\right) \oplus (F_{-1} \cap H^*P_{-1})$. By construction $H^*L_0 \cong F_{-1} \cap H^*P_{-1}$, hence Proposition~\ref{prop:quotientSA2} yields the result. 	
\end{proof}

\begin{prop} \label{prop:realizeinclusionL0-P-1}
	There is a $\tmf$-module map $\iota: \tmf \sm P_{-1} \to \tmf \sm L_0$ realizing the inclusion $\overline{\iota} : H^* L_0 \to H^* P_{-1}$ of $\SA(2)$-modules. Explicitly, 
	\begin{equation*}
	H^*\iota : \SA\sur\SA(2) \otimes H^* L_0 \to \SA\sur\SA(2) \otimes H^* P_{-1}
	\end{equation*}
	is $\SA\sur\SA(2) \otimes \overline{\iota}$.
\end{prop}	

\begin{proof}
	By adjunction, it suffices to show that the corresponding map of spectra $P_{-1} \to \tmf \sm L_0$ survives the Adams spectral sequence.  By the change-of-rings isomorphism \eqref{eq:changeofrings} the Adams $E_2$-page computing the desired homotopy classes is
	\begin{equation}\label{eq:E2P-1toL0}
		E_2 \cong \Ext_{\SA(2)}^{s,t}(H^*L_0, x^{-1} \FF_2[x])
	\end{equation}
	One concludes the argument by observing that the relevant extension group vanishes whenever $t-s = -1$. Indeed, there is a spectral sequence
	\begin{equation}\label{eq:E1P-1toL0}
	E_1 \cong \Ext_{\SA(2)}^{s,t}(H^*L_0, H^*L_0) \oplus \bigoplus_{j \geq 0}\Ext_{\SA(2)}^{s,t}(H^*L_0, \sus{8j-1}\SA(2)\sur\SA(1))
	\end{equation}	
	which abuts to the $E_2$-term \eqref{eq:E2P-1toL0}. It then suffices to check there is nothing in degree $t-s = -1$ in \eqref{eq:E1P-1toL0}.
	
	The first summand, $\Ext_{\SA(2)}^{s,t}(H^* L_0, H^* L_0)$, requires a direct computation. The relevant $\Ext$ chart is displayed in Figure~\ref{fig:endl0extchart} clearly has no classes in stem $t - s = -1$.
	
	The observant reader will note that, as an $\SA(2)$-module, the vector space dual of $\SA(2)\sur\SA(1)$ is $\sus{-17} \SA(2)\sur\SA(1)$. In particular, by adjunction and change-of-rings, the second summand is isomorphic to $\bigoplus_{j\geq 0} \Ext_{\SA(1)}^{s,t}(\sus{-8(j+2)} H^*L_0, \FF_2)$. This is a straightforward computation in $\SA(1)$-modules, since $H^* L_0 = \FF_2 \oplus \Sigma QM \oplus \sus{8} \FF_2$ where $QM$ is the question mark complex, yielding the $E_2$-term of $\ko \vee \ko\langle1\rangle \vee \sus{8}\ko$ which are all 8-periodic.  Figure~\ref{fig:ExtA1L0} displays $\Ext_{\SA(1)}^{s,t}(H^* L_0, \FF_2)$, and shows there cannot be a class in stem $t-s=-1$ for degree reasons.
\end{proof}

\begin{defn} \label{defn:ttmf-1}
	Choose a $\tmf$-module map $\pi$ satisfying the hypothesis of Proposition~\ref{prop:realizeinclusionL0-P-1}. Let $t(\tmf)_{-1}$ be the fiber of $\pi$.
\end{defn}

\section{The Tate spectrum of $\tmf$} \label{sec:Tatespectrum}

The proof of Theorem~\ref{thm:maintheorem} is decomposed into a series of lemmas. The argument presented here is similar to that given for the splitting of $t(\ko)$ \cite{DM84}.  In what follows, we denote $\Sigma t(\tmf)_{-1}$ by $\overline{t(\tmf)}$. We will show that $\overline{t(\tmf)}$ splits as a wedge of copies of $\ko$.  This is done in two steps: first, we show that the cohomology of $\overline{t(\tmf)}$ splits as a module over $\SA$ in Lemma~\ref{lem:ttmfbarsplits}; second, we compute the $E_2$-page of the Adams spectral sequence converging to $[\overline{t(\tmf)},\ko]$ to show that there are classes in $[\overline{t(\tmf)},\sus{8j}\ko]$ for all $j \in \ZZ$ that realize the splitting (this is done in Lemma~\ref{lem:ttmfko}).  The main result we need is the following computation of Mahowald~\cite{Mah81}:

\begin{lem}\label{lem:stablesplit}
	The stable $\SA(1)$-module $\SA\sur\SA(1)$ splits (in the stable category of $\SA(1)$-modules) as
	\[
		\bigoplus_{\ell \geq 0} \sus{12\ell+\sigma(\ell)}\dsus{4\ell - \sigma(\ell)}\FF_2 \oplus \bigoplus_{\ell\geq 0} \sus{12\ell+\sigma(\ell)+4}\dsus{4\ell - \sigma(\ell)}\Lambda(Sq^2).
	\]
\end{lem}

Note that, in particular, this gives that
\begin{align*}
\Ext_{\SA(1)}^{s,t}(\FF_2,\SA\sur\SA(1)) \cong &\bigoplus_{\ell \geq 0} \Ext_{\SA(1)}^{s+4\ell-\sigma(\ell),t+12\ell-\sigma(\ell)}(\FF_2,\FF_2) \oplus \\ &\bigoplus_{\ell \geq 0} \Ext_{E(1)}^{s+4\ell-\sigma(\ell),t+12\ell-\sigma(\ell)+4}(\FF_2,\FF_2)
\end{align*}
for $s > 0$, and these groups are zero as soon as $t - s \equiv 3 \mod 4$.

\begin{lem} \label{lem:ttmfbarsplits}
	There is an isomorphism of $\SA$-modules $H^*\overline{t(\tmf)} \cong \SA\sur\SA(1)[x^8]$ where $x^8$ is in degree 8.
\end{lem}

\begin{proof} 
	By Proposition~\ref{prop:associategradedP-1}, Proposition~\ref{prop:realizeinclusionL0-P-1}, and Definition~\ref{defn:ttmf-1} there is a filtration of $H^*\overline{t(\tmf)}$ whose associated graded is $\bigoplus_{j\geq 0} \sus{8j-1} \SA(2)\sur\SA(1)$. To conclude, we show inductively that there are no non-trivial extensions
	\[
		\Ext_\SA^{1,0}\left( \bigoplus_{j=0}^n \sus{8j} \SA\sur\SA(1), \SA\sur\SA(1) \right) = \Ext_{\SA(1)}^{1,t+8}(\FF_2, \SA\sur\SA(1)[x^8]/(x^n))
	\]
	By Lemma~\ref{lem:stablesplit}, this $\Ext$ vanishes whenever $t-s\equiv -1 \mod 4$. The degrees in which we are looking for non-trivial extensions are of the form $(1,8j)$. The result follows.
\end{proof}

Note that, in particular, the Adams spectral sequence computing $\pi_*(\overline{t(\tmf)})$ gives homotopy classes $x^{8j} \in \pi_{8j}(\overline{t(\tmf)}$. Indeed, the $E_2$-page of this spectral sequence is given by
\begin{equation}\label{eq:E2ttmf}
E_2 = \Ext_\SA(H^*\overline{t(\tmf)}, \FF_2) \cong \Ext_{\SA(1)}(\FF_2[x^8],\FF_2)
\end{equation}
and Lemma~\ref{lem:stablesplit} ensures that the elements of $x^{8j}$ survive the Adams spectral sequence for degree reasons.

\begin{lem}\label{lem:commdiag}
	We can arrange the cofiber sequences
	\begin{equation}\label{eq:fiberseq}
	\sus{-8k} t(\tmf)_{-1} \to \tmf \sm P_{-1-8k} \to \sus{-8k}\tmf \sm L_0
	\end{equation}
	in a commutative diagram
	\begin{equation}\label{eq:commdiag}
	\begin{CD}
		\sus{-8k} \overline{t(\tmf)} @>>> \Sigma \tmf \sm P_{-1-8k} @>>> \sus{-8k+1}\tmf \sm L_0 \\
		@VVV @VVV @VVV\\
		\sus{-8k + 8} \overline{t(\tmf)} @>>> \Sigma \tmf \sm P_{7-8k} @>>> \sus{-8k+9}\tmf \sm L_0
	\end{CD}
	\end{equation}
	Moreover, $t(\tmf)$ is the limit of both the leftmost and middle terms.	
\end{lem}

\begin{proof}
	The existence of the left commutative square in \eqref{eq:commdiag} comes from the identification $\sus{8} \tmf \sm P_n \simeq \tmf \sm P_{n+8}$ for all $n \in \ZZ$ \cite{BDM00}. This gives the asserted compatibility between the cofiber sequences.
	
	We now take the inverse limit on the three terms. By definition, the limit of the middle term is $\sus{-1}t(\tmf)$. We need to show that the map $\varprojlim \sus{-1}\tmf\sm P_{-1-8k} \to t(\tmf)$ is a weak equivalence. Note that the composite of the right-most vertical maps in \eqref{eq:commdiag} is zero, since it belongs to $[L_0, \tmf \sm L_0]_{-16} = 0$, by adjunction. Thus the limit is contractible.
\end{proof}

In light of Lemma~\ref{lem:commdiag}, we reduce our analysis of $t(\tmf)$ to $\overline{t(\tmf)}$. This, in turn, reduces to a simple Adams spectral sequence computation.

\begin{lem} \label{lem:ttmfko}
	There is a weak equivalence $\overline{t(\tmf)} \simeq \ko[x^8]$. Moreover, the maps $\sus{-8k} \overline{t(\tmf)} \to \sus{-8k+8} \overline{t(\tmf)}$ coincide with multiplication by $x^8$.
\end{lem}

\begin{proof}
We argue as follows: we build a map $\phi_j : \overline{t(\tmf)} \to \sus{8j} \ko$ for all $j \geq 0$, which realizes the injection $\phi^*_j : \sus{8j}\SA\sur\SA(1) \to \SA\sur\SA(1)[x^8]$. Then the coproduct of the $\phi_j$ will be the desired weak equivalence.
	
	First, fix $j \geq 0$ and build the map $\phi_j$. To this end, we compute the Adams spectral sequence converging to $[\overline{t(\tmf)}, \sus{8j}\ko]$. Its $E_2$-page is $\Ext_\SA^{s,t+8j}(\SA\sur\SA(1),\SA\sur\SA(1))[x^{\pm8}]$. In particular, Lemma~\ref{lem:stablesplit} guarantees an element in $(s,t+8j) = (0,8j)$ cannot be the source of a differential for degree reasons. This precisely means that any $\SA$-module map $\sus{8j}H^* \ko \to H^*\overline{t(\tmf)}$ comes from a map $\overline{t(\tmf)} \to \sus{8j}\ko$. Choose one map corresponding to the inclusion of $\sus{8j}\SA\sur\SA(1)$ into $\SA\sur\SA(1)[x^{\pm 8}]$ and call it $\phi_j$.
	
	Let $\phi: \overline{t(\tmf)} \to \ko[x^8]$ be the wedge of the above $\phi_j$ for $j \geq 0$. By construction, it induces an isomorphism in cohomology between connective spectra, and thus it is a (2-local) equivalence.
	
	The assertion about the maps $\sus{-8k} \overline{t(\tmf)} \to \sus{-8k+8}\overline{t(\tmf)}$ follows from its compatibility with the decomposition $\Ext_\SA(\SA\sur\SA(1)\otimes \SA\sur\SA(1),\FF_2)[x^8]$. 
\end{proof}

\begin{proof}[Proof of Theorem~\ref{thm:maintheorem}]
	By Lemma~\ref{lem:commdiag}, $t(\tmf) \simeq \varprojlim \sus{-8k}\overline{t(\tmf)}$ which, by Lemma~\ref{lem:ttmfko}, is weakly equivalent to $\varprojlim \sus{-8k} \ko[x^8] \simeq \ko[x^{\pm 8}]$.
\end{proof}

\begin{figure}
	
	\begin{tikzpicture}[scale=0.6]
	
	\clip (-1,-1) rectangle ( 20.50, 10.50);
	\draw[color=lightgray] (0,0) grid [step=4]  (20,10);
	
	\node [below] at (0,0) {$-9$};
	\node [below] at (4,0) {$-5$};
	\node [below] at (8,0) {$-1$};
	\node [below] at (12,0) {$3$};
	\node [below] at (16,0) {$7$};
	\node [below] at (20,0) {$11$};
	
	\node [left] at (0,0) {$0$};
	\node [left] at (0,4) {$4$};
	\node [left] at (0,8) {$8$};
	
	\draw [fill] ( 1.00, 0.00) circle [radius=0.05];
	\draw [fill] ( 2.00, 0.00) circle [radius=0.05];
	\draw [fill] ( 9.07,-0.07) circle [radius=0.05];
	\draw [fill] ( 8.93, 0.07) circle [radius=0.05];
	\draw [fill] (10.00, 0.00) circle [radius=0.05];
	\draw [fill] ( 1.00, 1.00) circle [radius=0.05];
	\draw [fill] ( 2.00, 1.00) circle [radius=0.05];
	\draw [fill] ( 3.00, 1.00) circle [radius=0.05];
	\draw [fill] ( 9.07, 0.93) circle [radius=0.05];
	\draw [fill] ( 8.93, 1.07) circle [radius=0.05];
	\draw [fill] (10.14, 0.86) circle [radius=0.05];
	\draw [fill] (10.00, 1.00) circle [radius=0.05];
	\draw [fill] ( 9.86, 1.14) circle [radius=0.05];
	\draw [fill] (11.07, 0.93) circle [radius=0.05];
	\draw [fill] (10.93, 1.07) circle [radius=0.05];
	\draw [fill] (12.07, 0.93) circle [radius=0.05];
	\draw [fill] (11.93, 1.07) circle [radius=0.05];
	\draw [fill] (13.00, 1.00) circle [radius=0.05];
	\draw [fill] ( 1.00, 2.00) circle [radius=0.05];
	\draw [fill] ( 3.00, 2.00) circle [radius=0.05];
	\draw [fill] ( 5.00, 2.00) circle [radius=0.05];
	\draw [fill] ( 9.07, 1.93) circle [radius=0.05];
	\draw [fill] ( 8.93, 2.07) circle [radius=0.05];
	\draw [fill] (11.14, 1.86) circle [radius=0.05];
	\draw [fill] (11.00, 2.00) circle [radius=0.05];
	\draw [fill] (10.86, 2.14) circle [radius=0.05];
	\draw [fill] (12.07, 1.93) circle [radius=0.05];
	\draw [fill] (11.93, 2.07) circle [radius=0.05];
	\draw [fill] (13.07, 1.93) circle [radius=0.05];
	\draw [fill] (12.93, 2.07) circle [radius=0.05];
	\draw [fill] (14.00, 2.00) circle [radius=0.05];
	\draw [fill] (15.07, 1.93) circle [radius=0.05];
	\draw [fill] (14.93, 2.07) circle [radius=0.05];
	\draw [fill] (16.00, 2.00) circle [radius=0.05];
	\draw [fill] ( 1.00, 3.00) circle [radius=0.05];
	\draw [fill] ( 5.07, 2.93) circle [radius=0.05];
	\draw [fill] ( 4.93, 3.07) circle [radius=0.05];
	\draw [fill] ( 9.14, 2.86) circle [radius=0.05];
	\draw [fill] ( 9.00, 3.00) circle [radius=0.05];
	\draw [fill] ( 8.86, 3.14) circle [radius=0.05];
	\draw [fill] (12.00, 3.00) circle [radius=0.05];
	\draw [fill] (13.14, 2.86) circle [radius=0.05];
	\draw [fill] (13.00, 3.00) circle [radius=0.05];
	\draw [fill] (12.86, 3.14) circle [radius=0.05];
	\draw [fill] (16.00, 3.00) circle [radius=0.05];
	\draw [fill] (17.21, 2.79) circle [radius=0.05];
	\draw [fill] (17.07, 2.93) circle [radius=0.05];
	\draw [fill] (16.93, 3.07) circle [radius=0.05];
	\draw [fill] (16.79, 3.21) circle [radius=0.05];
	\draw [fill] ( 1.00, 4.00) circle [radius=0.05];
	\draw [fill] ( 5.07, 3.93) circle [radius=0.05];
	\draw [fill] ( 4.93, 4.07) circle [radius=0.05];
	\draw [fill] ( 9.21, 3.79) circle [radius=0.05];
	\draw [fill] ( 9.07, 3.93) circle [radius=0.05];
	\draw [fill] ( 8.93, 4.07) circle [radius=0.05];
	\draw [fill] ( 8.79, 4.21) circle [radius=0.05];
	\draw [fill] (10.00, 4.00) circle [radius=0.05];
	\draw [fill] (13.14, 3.86) circle [radius=0.05];
	\draw [fill] (13.00, 4.00) circle [radius=0.05];
	\draw [fill] (12.86, 4.14) circle [radius=0.05];
	\draw [fill] (17.28, 3.72) circle [radius=0.05];
	\draw [fill] (17.14, 3.86) circle [radius=0.05];
	\draw [fill] (17.00, 4.00) circle [radius=0.05];
	\draw [fill] (16.86, 4.14) circle [radius=0.05];
	\draw [fill] (16.72, 4.28) circle [radius=0.05];
	\draw [fill] (18.14, 3.86) circle [radius=0.05];
	\draw [fill] (18.00, 4.00) circle [radius=0.05];
	\draw [fill] (17.86, 4.14) circle [radius=0.05];
	\draw [fill] ( 1.00, 5.00) circle [radius=0.05];
	\draw [fill] ( 5.07, 4.93) circle [radius=0.05];
	\draw [fill] ( 4.93, 5.07) circle [radius=0.05];
	\draw [fill] ( 9.21, 4.79) circle [radius=0.05];
	\draw [fill] ( 9.07, 4.93) circle [radius=0.05];
	\draw [fill] ( 8.93, 5.07) circle [radius=0.05];
	\draw [fill] ( 8.79, 5.21) circle [radius=0.05];
	\draw [fill] (10.00, 5.00) circle [radius=0.05];
	\draw [fill] (11.00, 5.00) circle [radius=0.05];
	\draw [fill] (13.14, 4.86) circle [radius=0.05];
	\draw [fill] (13.00, 5.00) circle [radius=0.05];
	\draw [fill] (12.86, 5.14) circle [radius=0.05];
	\draw [fill] (17.28, 4.72) circle [radius=0.05];
	\draw [fill] (17.14, 4.86) circle [radius=0.05];
	\draw [fill] (17.00, 5.00) circle [radius=0.05];
	\draw [fill] (16.86, 5.14) circle [radius=0.05];
	\draw [fill] (16.72, 5.28) circle [radius=0.05];
	\draw [fill] (18.14, 4.86) circle [radius=0.05];
	\draw [fill] (18.00, 5.00) circle [radius=0.05];
	\draw [fill] (17.86, 5.14) circle [radius=0.05];
	\draw [fill] (19.07, 4.93) circle [radius=0.05];
	\draw [fill] (18.93, 5.07) circle [radius=0.05];
	\draw [fill] (20.00, 5.00) circle [radius=0.05];
	\draw [fill] ( 1.00, 6.00) circle [radius=0.05];
	\draw [fill] ( 5.07, 5.93) circle [radius=0.05];
	\draw [fill] ( 4.93, 6.07) circle [radius=0.05];
	\draw [fill] ( 9.21, 5.79) circle [radius=0.05];
	\draw [fill] ( 9.07, 5.93) circle [radius=0.05];
	\draw [fill] ( 8.93, 6.07) circle [radius=0.05];
	\draw [fill] ( 8.79, 6.21) circle [radius=0.05];
	\draw [fill] (11.00, 6.00) circle [radius=0.05];
	\draw [fill] (13.21, 5.79) circle [radius=0.05];
	\draw [fill] (13.07, 5.93) circle [radius=0.05];
	\draw [fill] (12.93, 6.07) circle [radius=0.05];
	\draw [fill] (12.79, 6.21) circle [radius=0.05];
	\draw [fill] (17.28, 5.72) circle [radius=0.05];
	\draw [fill] (17.14, 5.86) circle [radius=0.05];
	\draw [fill] (17.00, 6.00) circle [radius=0.05];
	\draw [fill] (16.86, 6.14) circle [radius=0.05];
	\draw [fill] (16.72, 6.28) circle [radius=0.05];
	\draw [fill] (19.14, 5.86) circle [radius=0.05];
	\draw [fill] (19.00, 6.00) circle [radius=0.05];
	\draw [fill] (18.86, 6.14) circle [radius=0.05];
	\draw [fill] (20.00, 6.00) circle [radius=0.05];
	\draw [fill] ( 1.00, 7.00) circle [radius=0.05];
	\draw [fill] ( 5.07, 6.93) circle [radius=0.05];
	\draw [fill] ( 4.93, 7.07) circle [radius=0.05];
	\draw [fill] ( 9.21, 6.79) circle [radius=0.05];
	\draw [fill] ( 9.07, 6.93) circle [radius=0.05];
	\draw [fill] ( 8.93, 7.07) circle [radius=0.05];
	\draw [fill] ( 8.79, 7.21) circle [radius=0.05];
	\draw [fill] (13.28, 6.72) circle [radius=0.05];
	\draw [fill] (13.14, 6.86) circle [radius=0.05];
	\draw [fill] (13.00, 7.00) circle [radius=0.05];
	\draw [fill] (12.86, 7.14) circle [radius=0.05];
	\draw [fill] (12.72, 7.28) circle [radius=0.05];
	\draw [fill] (17.35, 6.65) circle [radius=0.05];
	\draw [fill] (17.21, 6.79) circle [radius=0.05];
	\draw [fill] (17.07, 6.93) circle [radius=0.05];
	\draw [fill] (16.93, 7.07) circle [radius=0.05];
	\draw [fill] (16.79, 7.21) circle [radius=0.05];
	\draw [fill] (16.65, 7.35) circle [radius=0.05];
	\draw [fill] (20.00, 7.00) circle [radius=0.05];
	\draw [fill] ( 1.00, 8.00) circle [radius=0.05];
	\draw [fill] ( 5.07, 7.93) circle [radius=0.05];
	\draw [fill] ( 4.93, 8.07) circle [radius=0.05];
	\draw [fill] ( 9.21, 7.79) circle [radius=0.05];
	\draw [fill] ( 9.07, 7.93) circle [radius=0.05];
	\draw [fill] ( 8.93, 8.07) circle [radius=0.05];
	\draw [fill] ( 8.79, 8.21) circle [radius=0.05];
	\draw [fill] (13.28, 7.72) circle [radius=0.05];
	\draw [fill] (13.14, 7.86) circle [radius=0.05];
	\draw [fill] (13.00, 8.00) circle [radius=0.05];
	\draw [fill] (12.86, 8.14) circle [radius=0.05];
	\draw [fill] (12.72, 8.28) circle [radius=0.05];
	\draw [fill] (17.42, 7.58) circle [radius=0.05];
	\draw [fill] (17.28, 7.72) circle [radius=0.05];
	\draw [fill] (17.14, 7.86) circle [radius=0.05];
	\draw [fill] (17.00, 8.00) circle [radius=0.05];
	\draw [fill] (16.86, 8.14) circle [radius=0.05];
	\draw [fill] (16.72, 8.28) circle [radius=0.05];
	\draw [fill] (16.58, 8.42) circle [radius=0.05];
	\draw [fill] (18.00, 8.00) circle [radius=0.05];
	\draw [fill] ( 1.00, 9.00) circle [radius=0.05];
	\draw [fill] ( 5.07, 8.93) circle [radius=0.05];
	\draw [fill] ( 4.93, 9.07) circle [radius=0.05];
	\draw [fill] ( 9.21, 8.79) circle [radius=0.05];
	\draw [fill] ( 9.07, 8.93) circle [radius=0.05];
	\draw [fill] ( 8.93, 9.07) circle [radius=0.05];
	\draw [fill] ( 8.79, 9.21) circle [radius=0.05];
	\draw [fill] (13.28, 8.72) circle [radius=0.05];
	\draw [fill] (13.14, 8.86) circle [radius=0.05];
	\draw [fill] (13.00, 9.00) circle [radius=0.05];
	\draw [fill] (12.86, 9.14) circle [radius=0.05];
	\draw [fill] (12.72, 9.28) circle [radius=0.05];
	\draw [fill] (17.42, 8.58) circle [radius=0.05];
	\draw [fill] (17.28, 8.72) circle [radius=0.05];
	\draw [fill] (17.14, 8.86) circle [radius=0.05];
	\draw [fill] (17.00, 9.00) circle [radius=0.05];
	\draw [fill] (16.86, 9.14) circle [radius=0.05];
	\draw [fill] (16.72, 9.28) circle [radius=0.05];
	\draw [fill] (16.58, 9.42) circle [radius=0.05];
	\draw [fill] (18.00, 9.00) circle [radius=0.05];
	\draw [fill] (19.00, 9.00) circle [radius=0.05];
	\draw [fill] ( 1.00,10.00) circle [radius=0.05];
	\draw [fill] ( 5.07, 9.93) circle [radius=0.05];
	\draw [fill] ( 4.93,10.07) circle [radius=0.05];
	\draw [fill] ( 9.21, 9.79) circle [radius=0.05];
	\draw [fill] ( 9.07, 9.93) circle [radius=0.05];
	\draw [fill] ( 8.93,10.07) circle [radius=0.05];
	\draw [fill] ( 8.79,10.21) circle [radius=0.05];
	\draw [fill] (13.28, 9.72) circle [radius=0.05];
	\draw [fill] (13.14, 9.86) circle [radius=0.05];
	\draw [fill] (13.00,10.00) circle [radius=0.05];
	\draw [fill] (12.86,10.14) circle [radius=0.05];
	\draw [fill] (12.72,10.28) circle [radius=0.05];
	\draw [fill] (17.42, 9.58) circle [radius=0.05];
	\draw [fill] (17.28, 9.72) circle [radius=0.05];
	\draw [fill] (17.14, 9.86) circle [radius=0.05];
	\draw [fill] (17.00,10.00) circle [radius=0.05];
	\draw [fill] (16.86,10.14) circle [radius=0.05];
	\draw [fill] (16.72,10.28) circle [radius=0.05];
	\draw [fill] (16.58,10.42) circle [radius=0.05];
	\draw [fill] (19.00,10.00) circle [radius=0.05];

	\draw ( 1.00, 1.00) --( 1.00, 0.00);
	\draw ( 2.00, 1.00) --( 1.00, 0.00);
	\draw ( 3.00, 1.00) --( 2.00, 0.00);
	\draw ( 9.07, 0.93) --( 9.07,-0.07);
	\draw ( 8.93, 1.07) --( 8.93, 0.07);
	\draw (10.00, 1.00) --( 9.07,-0.07);
	\draw ( 9.86, 1.14) --( 8.93, 0.07);
	\draw (10.93, 1.07) --(10.00, 0.00);
	\draw [dashed]  (12.07, 0.93) --( 9.07,-0.07);
	\draw [dashed]  (11.93, 1.07) --( 8.93, 0.07);
	\draw [dashed]  (13.00, 1.00) --(10.00, 0.00);
	\draw ( 1.00, 2.00) --( 1.00, 1.00);
	\draw ( 3.00, 2.00) --( 2.00, 1.00);
	\draw ( 9.07, 1.93) --( 9.07, 0.93);
	\draw ( 8.93, 2.07) --( 8.93, 1.07);
	\draw (11.14, 1.86) --(10.14, 0.86);
	\draw (11.00, 2.00) --(10.00, 1.00);
	\draw (10.86, 2.14) --( 9.86, 1.14);
	\draw [dashed]  (12.07, 1.93) --( 9.07, 0.93);
	\draw (12.07, 1.93) --(12.07, 0.93);
	\draw [dashed]  (11.93, 2.07) --( 8.93, 1.07);
	\draw (11.93, 2.07) --(11.07, 0.93);
	\draw (11.93, 2.07) --(11.93, 1.07);
	\draw [dashed]  (13.07, 1.93) --(10.14, 0.86);
	\draw [dashed]  (14.00, 2.00) --(11.07, 0.93);
	\draw [dashed]  (15.07, 1.93) --(12.07, 0.93);
	\draw [dashed]  (14.93, 2.07) --(11.93, 1.07);
	\draw [dashed]  (16.00, 2.00) --(13.00, 1.00);
	\draw ( 1.00, 3.00) --( 1.00, 2.00);
	\draw ( 4.93, 3.07) --( 5.00, 2.00);
	\draw ( 9.00, 3.00) --( 9.07, 1.93);
	\draw ( 8.86, 3.14) --( 8.93, 2.07);
	\draw [dashed]  (12.00, 3.00) --( 9.07, 1.93);
	\draw (12.00, 3.00) --(11.00, 2.00);
	\draw (12.00, 3.00) --(12.07, 1.93);
	\draw (12.86, 3.14) --(12.93, 2.07);
	\draw [dashed]  (16.00, 3.00) --(13.07, 1.93);
	\draw [dashed]  (16.79, 3.21) --(14.00, 2.00);
	\draw ( 1.00, 4.00) --( 1.00, 3.00);
	\draw ( 5.07, 3.93) --( 5.07, 2.93);
	\draw ( 4.93, 4.07) --( 4.93, 3.07);
	\draw ( 9.07, 3.93) --( 9.14, 2.86);
	\draw ( 8.93, 4.07) --( 9.00, 3.00);
	\draw ( 8.79, 4.21) --( 8.86, 3.14);
	\draw (10.00, 4.00) --( 9.14, 2.86);
	\draw (13.14, 3.86) --(13.14, 2.86);
	\draw (13.00, 4.00) --(13.00, 3.00);
	\draw (12.86, 4.14) --(12.86, 3.14);
	\draw (16.86, 4.14) --(17.21, 2.79);
	\draw (16.72, 4.28) --(17.07, 2.93);
	\draw (18.14, 3.86) --(17.21, 2.79);
	\draw (18.00, 4.00) --(17.07, 2.93);
	\draw (17.86, 4.14) --(16.93, 3.07);
	\draw ( 1.00, 5.00) --( 1.00, 4.00);
	\draw ( 5.07, 4.93) --( 5.07, 3.93);
	\draw ( 4.93, 5.07) --( 4.93, 4.07);
	\draw ( 9.21, 4.79) --( 9.21, 3.79);
	\draw ( 9.07, 4.93) --( 9.07, 3.93);
	\draw ( 8.93, 5.07) --( 8.93, 4.07);
	\draw ( 8.79, 5.21) --( 8.79, 4.21);
	\draw (10.00, 5.00) --( 9.21, 3.79);
	\draw (11.00, 5.00) --(10.00, 4.00);
	\draw (13.14, 4.86) --(13.14, 3.86);
	\draw (13.00, 5.00) --(13.00, 4.00);
	\draw (12.86, 5.14) --(12.86, 4.14);
	\draw (17.28, 4.72) --(17.28, 3.72);
	\draw (17.14, 4.86) --(17.14, 3.86);
	\draw (17.00, 5.00) --(17.00, 4.00);
	\draw (16.86, 5.14) --(16.86, 4.14);
	\draw (16.72, 5.28) --(16.72, 4.28);
	\draw (18.14, 4.86) --(17.28, 3.72);
	\draw (18.00, 5.00) --(17.14, 3.86);
	\draw (17.86, 5.14) --(17.00, 4.00);
	\draw (19.07, 4.93) --(18.14, 3.86);
	\draw (18.93, 5.07) --(18.00, 4.00);
	\draw [dashed]  (20.00, 5.00) --(17.28, 3.72);
	\draw ( 1.00, 6.00) --( 1.00, 5.00);
	\draw ( 5.07, 5.93) --( 5.07, 4.93);
	\draw ( 4.93, 6.07) --( 4.93, 5.07);
	\draw ( 9.21, 5.79) --( 9.21, 4.79);
	\draw ( 9.07, 5.93) --( 9.07, 4.93);
	\draw ( 8.93, 6.07) --( 8.93, 5.07);
	\draw ( 8.79, 6.21) --( 8.79, 5.21);
	\draw (11.00, 6.00) --(10.00, 5.00);
	\draw (13.07, 5.93) --(13.14, 4.86);
	\draw (12.93, 6.07) --(13.00, 5.00);
	\draw (12.79, 6.21) --(12.86, 5.14);
	\draw (17.28, 5.72) --(17.28, 4.72);
	\draw (17.14, 5.86) --(17.14, 4.86);
	\draw (17.00, 6.00) --(17.00, 5.00);
	\draw (16.86, 6.14) --(16.86, 5.14);
	\draw (16.72, 6.28) --(16.72, 5.28);
	\draw (19.14, 5.86) --(18.14, 4.86);
	\draw (19.00, 6.00) --(18.00, 5.00);
	\draw (18.86, 6.14) --(17.86, 5.14);
	\draw [dashed]  (20.00, 6.00) --(17.28, 4.72);
	\draw (20.00, 6.00) --(20.00, 5.00);
	\draw [dashed]  (23.00, 6.00) --(20.00, 5.00);
	\draw ( 1.00, 7.00) --( 1.00, 6.00);
	\draw ( 5.07, 6.93) --( 5.07, 5.93);
	\draw ( 4.93, 7.07) --( 4.93, 6.07);
	\draw ( 9.21, 6.79) --( 9.21, 5.79);
	\draw ( 9.07, 6.93) --( 9.07, 5.93);
	\draw ( 8.93, 7.07) --( 8.93, 6.07);
	\draw ( 8.79, 7.21) --( 8.79, 6.21);
	\draw (13.14, 6.86) --(13.21, 5.79);
	\draw (13.00, 7.00) --(13.07, 5.93);
	\draw (12.86, 7.14) --(12.93, 6.07);
	\draw (12.72, 7.28) --(12.79, 6.21);
	\draw (17.21, 6.79) --(17.28, 5.72);
	\draw (17.07, 6.93) --(17.14, 5.86);
	\draw (16.93, 7.07) --(17.00, 6.00);
	\draw (16.79, 7.21) --(16.86, 6.14);
	\draw (16.65, 7.35) --(16.72, 6.28);
	\draw [dashed]  (20.00, 7.00) --(17.28, 5.72);
	\draw (20.00, 7.00) --(19.14, 5.86);
	\draw (20.00, 7.00) --(20.00, 6.00);
	\draw ( 1.00, 8.00) --( 1.00, 7.00);
	\draw ( 5.07, 7.93) --( 5.07, 6.93);
	\draw ( 4.93, 8.07) --( 4.93, 7.07);
	\draw ( 9.21, 7.79) --( 9.21, 6.79);
	\draw ( 9.07, 7.93) --( 9.07, 6.93);
	\draw ( 8.93, 8.07) --( 8.93, 7.07);
	\draw ( 8.79, 8.21) --( 8.79, 7.21);
	\draw (13.28, 7.72) --(13.28, 6.72);
	\draw (13.14, 7.86) --(13.14, 6.86);
	\draw (13.00, 8.00) --(13.00, 7.00);
	\draw (12.86, 8.14) --(12.86, 7.14);
	\draw (12.72, 8.28) --(12.72, 7.28);
	\draw (17.28, 7.72) --(17.35, 6.65);
	\draw (17.14, 7.86) --(17.21, 6.79);
	\draw (17.00, 8.00) --(17.07, 6.93);
	\draw (16.86, 8.14) --(16.93, 7.07);
	\draw (16.72, 8.28) --(16.79, 7.21);
	\draw (16.58, 8.42) --(16.65, 7.35);
	\draw (18.00, 8.00) --(17.35, 6.65);
	\draw ( 1.00, 9.00) --( 1.00, 8.00);
	\draw ( 5.07, 8.93) --( 5.07, 7.93);
	\draw ( 4.93, 9.07) --( 4.93, 8.07);
	\draw ( 9.21, 8.79) --( 9.21, 7.79);
	\draw ( 9.07, 8.93) --( 9.07, 7.93);
	\draw ( 8.93, 9.07) --( 8.93, 8.07);
	\draw ( 8.79, 9.21) --( 8.79, 8.21);
	\draw (13.28, 8.72) --(13.28, 7.72);
	\draw (13.14, 8.86) --(13.14, 7.86);
	\draw (13.00, 9.00) --(13.00, 8.00);
	\draw (12.86, 9.14) --(12.86, 8.14);
	\draw (12.72, 9.28) --(12.72, 8.28);
	\draw (17.42, 8.58) --(17.42, 7.58);
	\draw (17.28, 8.72) --(17.28, 7.72);
	\draw (17.14, 8.86) --(17.14, 7.86);
	\draw (17.00, 9.00) --(17.00, 8.00);
	\draw (16.86, 9.14) --(16.86, 8.14);
	\draw (16.72, 9.28) --(16.72, 8.28);
	\draw (16.58, 9.42) --(16.58, 8.42);
	\draw (18.00, 9.00) --(17.42, 7.58);
	\draw (19.00, 9.00) --(18.00, 8.00);
	\draw ( 1.00,10.00) --( 1.00, 9.00);
	\draw ( 5.07, 9.93) --( 5.07, 8.93);
	\draw ( 4.93,10.07) --( 4.93, 9.07);
	\draw ( 9.21, 9.79) --( 9.21, 8.79);
	\draw ( 9.07, 9.93) --( 9.07, 8.93);
	\draw ( 8.93,10.07) --( 8.93, 9.07);
	\draw ( 8.79,10.21) --( 8.79, 9.21);
	\draw (13.28, 9.72) --(13.28, 8.72);
	\draw (13.14, 9.86) --(13.14, 8.86);
	\draw (13.00,10.00) --(13.00, 9.00);
	\draw (12.86,10.14) --(12.86, 9.14);
	\draw (12.72,10.28) --(12.72, 9.28);
	\draw (17.42, 9.58) --(17.42, 8.58);
	\draw (17.28, 9.72) --(17.28, 8.72);
	\draw (17.14, 9.86) --(17.14, 8.86);
	\draw (17.00,10.00) --(17.00, 9.00);
	\draw (16.86,10.14) --(16.86, 9.14);
	\draw (16.72,10.28) --(16.72, 9.28);
	\draw (16.58,10.42) --(16.58, 9.42);
	\draw (19.00,10.00) --(18.00, 9.00);
	\draw ( 1.00,11.00) --( 1.00,10.00);
	\draw ( 5.07,10.93) --( 5.07, 9.93);
	\draw ( 4.93,11.07) --( 4.93,10.07);
	\draw ( 9.21,10.79) --( 9.21, 9.79);
	\draw ( 9.07,10.93) --( 9.07, 9.93);
	\draw ( 8.93,11.07) --( 8.93,10.07);
	\draw ( 8.79,11.21) --( 8.79,10.21);
	\draw (13.28,10.72) --(13.28, 9.72);
	\draw (13.14,10.86) --(13.14, 9.86);
	\draw (13.00,11.00) --(13.00,10.00);
	\draw (12.86,11.14) --(12.86,10.14);
	\draw (12.72,11.28) --(12.72,10.28);
	\draw (17.42,10.58) --(17.42, 9.58);
	\draw (17.28,10.72) --(17.28, 9.72);
	\draw (17.14,10.86) --(17.14, 9.86);
	\draw (17.00,11.00) --(17.00,10.00);
	\draw (16.86,11.14) --(16.86,10.14);
	\draw (16.72,11.28) --(16.72,10.28);
	\draw (16.58,11.42) --(16.58,10.42);
	\end{tikzpicture}
	
	\caption{$\Ext^{s,t}_{\SA(2)}(H^*L_0, H^*L_0)$} \label{fig:endl0extchart}
\end{figure}

\begin{figure}
\begin{tikzpicture}[scale=0.6]

\clip (-1.5,-1) rectangle ( 21.50, 12.50);
\draw[color=lightgray] (0,0) grid [step=4]  (21,12);

\foreach \n in {-1,3,...,20}
{
	\def\nn{\n--1};
	\node [below] at (\nn,0) {$\n$};
}

\foreach \s in {0,4,...,12}
{
	\def\ss{\s-0};
	\node [left] at (-0.4,\ss,0) {$\s$};
}

\draw [fill] ( 1.00, 0.00) circle [radius=0.05];
\draw [fill] ( 2.00, 0.00) circle [radius=0.05];
\draw [fill] ( 9.00, 0.00) circle [radius=0.05];
\draw [fill] ( 1.00, 1.00) circle [radius=0.05];
\draw [fill] ( 2.00, 1.00) circle [radius=0.05];
\draw [fill] ( 3.00, 1.00) circle [radius=0.05];
\draw [fill] ( 9.00, 1.00) circle [radius=0.05];
\draw [fill] (10.00, 1.00) circle [radius=0.05];
\draw [fill] ( 1.00, 2.00) circle [radius=0.05];
\draw [fill] ( 3.00, 2.00) circle [radius=0.05];
\draw [fill] ( 5.00, 2.00) circle [radius=0.05];
\draw [fill] ( 9.00, 2.00) circle [radius=0.05];
\draw [fill] (11.00, 2.00) circle [radius=0.05];
\draw [fill] ( 1.00, 3.00) circle [radius=0.05];
\draw [fill] ( 5.07, 2.93) circle [radius=0.05];
\draw [fill] ( 4.93, 3.07) circle [radius=0.05];
\draw [fill] ( 9.07, 2.93) circle [radius=0.05];
\draw [fill] ( 8.93, 3.07) circle [radius=0.05];
\draw [fill] (13.00, 3.00) circle [radius=0.05];
\draw [fill] ( 1.00, 4.00) circle [radius=0.05];
\draw [fill] ( 5.07, 3.93) circle [radius=0.05];
\draw [fill] ( 4.93, 4.07) circle [radius=0.05];
\draw [fill] ( 9.14, 3.86) circle [radius=0.05];
\draw [fill] ( 9.00, 4.00) circle [radius=0.05];
\draw [fill] ( 8.86, 4.14) circle [radius=0.05];
\draw [fill] (10.00, 4.00) circle [radius=0.05];
\draw [fill] (13.00, 4.00) circle [radius=0.05];
\draw [fill] (17.00, 4.00) circle [radius=0.05];
\draw [fill] ( 1.00, 5.00) circle [radius=0.05];
\draw [fill] ( 5.07, 4.93) circle [radius=0.05];
\draw [fill] ( 4.93, 5.07) circle [radius=0.05];
\draw [fill] ( 9.14, 4.86) circle [radius=0.05];
\draw [fill] ( 9.00, 5.00) circle [radius=0.05];
\draw [fill] ( 8.86, 5.14) circle [radius=0.05];
\draw [fill] (10.00, 5.00) circle [radius=0.05];
\draw [fill] (11.00, 5.00) circle [radius=0.05];
\draw [fill] (13.00, 5.00) circle [radius=0.05];
\draw [fill] (17.00, 5.00) circle [radius=0.05];
\draw [fill] (18.00, 5.00) circle [radius=0.05];
\draw [fill] ( 1.00, 6.00) circle [radius=0.05];
\draw [fill] ( 5.07, 5.93) circle [radius=0.05];
\draw [fill] ( 4.93, 6.07) circle [radius=0.05];
\draw [fill] ( 9.14, 5.86) circle [radius=0.05];
\draw [fill] ( 9.00, 6.00) circle [radius=0.05];
\draw [fill] ( 8.86, 6.14) circle [radius=0.05];
\draw [fill] (11.00, 6.00) circle [radius=0.05];
\draw [fill] (13.07, 5.93) circle [radius=0.05];
\draw [fill] (12.93, 6.07) circle [radius=0.05];
\draw [fill] (17.00, 6.00) circle [radius=0.05];
\draw [fill] (19.00, 6.00) circle [radius=0.05];
\draw [fill] ( 1.00, 7.00) circle [radius=0.05];
\draw [fill] ( 5.07, 6.93) circle [radius=0.05];
\draw [fill] ( 4.93, 7.07) circle [radius=0.05];
\draw [fill] ( 9.14, 6.86) circle [radius=0.05];
\draw [fill] ( 9.00, 7.00) circle [radius=0.05];
\draw [fill] ( 8.86, 7.14) circle [radius=0.05];
\draw [fill] (13.14, 6.86) circle [radius=0.05];
\draw [fill] (13.00, 7.00) circle [radius=0.05];
\draw [fill] (12.86, 7.14) circle [radius=0.05];
\draw [fill] (17.07, 6.93) circle [radius=0.05];
\draw [fill] (16.93, 7.07) circle [radius=0.05];
\draw [fill] ( 1.00, 8.00) circle [radius=0.05];
\draw [fill] ( 5.07, 7.93) circle [radius=0.05];
\draw [fill] ( 4.93, 8.07) circle [radius=0.05];
\draw [fill] ( 9.14, 7.86) circle [radius=0.05];
\draw [fill] ( 9.00, 8.00) circle [radius=0.05];
\draw [fill] ( 8.86, 8.14) circle [radius=0.05];
\draw [fill] (13.14, 7.86) circle [radius=0.05];
\draw [fill] (13.00, 8.00) circle [radius=0.05];
\draw [fill] (12.86, 8.14) circle [radius=0.05];
\draw [fill] (17.14, 7.86) circle [radius=0.05];
\draw [fill] (17.00, 8.00) circle [radius=0.05];
\draw [fill] (16.86, 8.14) circle [radius=0.05];
\draw [fill] (18.00, 8.00) circle [radius=0.05];
\draw [fill] ( 1.00, 9.00) circle [radius=0.05];
\draw [fill] ( 5.07, 8.93) circle [radius=0.05];
\draw [fill] ( 4.93, 9.07) circle [radius=0.05];
\draw [fill] ( 9.14, 8.86) circle [radius=0.05];
\draw [fill] ( 9.00, 9.00) circle [radius=0.05];
\draw [fill] ( 8.86, 9.14) circle [radius=0.05];
\draw [fill] (13.14, 8.86) circle [radius=0.05];
\draw [fill] (13.00, 9.00) circle [radius=0.05];
\draw [fill] (12.86, 9.14) circle [radius=0.05];
\draw [fill] (17.14, 8.86) circle [radius=0.05];
\draw [fill] (17.00, 9.00) circle [radius=0.05];
\draw [fill] (16.86, 9.14) circle [radius=0.05];
\draw [fill] (18.00, 9.00) circle [radius=0.05];
\draw [fill] (19.00, 9.00) circle [radius=0.05];
\draw [fill] ( 1.00,10.00) circle [radius=0.05];
\draw [fill] ( 5.07, 9.93) circle [radius=0.05];
\draw [fill] ( 4.93,10.07) circle [radius=0.05];
\draw [fill] ( 9.14, 9.86) circle [radius=0.05];
\draw [fill] ( 9.00,10.00) circle [radius=0.05];
\draw [fill] ( 8.86,10.14) circle [radius=0.05];
\draw [fill] (13.14, 9.86) circle [radius=0.05];
\draw [fill] (13.00,10.00) circle [radius=0.05];
\draw [fill] (12.86,10.14) circle [radius=0.05];
\draw [fill] (17.14, 9.86) circle [radius=0.05];
\draw [fill] (17.00,10.00) circle [radius=0.05];
\draw [fill] (16.86,10.14) circle [radius=0.05];
\draw [fill] (19.00,10.00) circle [radius=0.05];
\draw [fill] ( 1.00,11.00) circle [radius=0.05];
\draw [fill] ( 5.07,10.93) circle [radius=0.05];
\draw [fill] ( 4.93,11.07) circle [radius=0.05];
\draw [fill] ( 9.14,10.86) circle [radius=0.05];
\draw [fill] ( 9.00,11.00) circle [radius=0.05];
\draw [fill] ( 8.86,11.14) circle [radius=0.05];
\draw [fill] (13.14,10.86) circle [radius=0.05];
\draw [fill] (13.00,11.00) circle [radius=0.05];
\draw [fill] (12.86,11.14) circle [radius=0.05];
\draw [fill] (17.14,10.86) circle [radius=0.05];
\draw [fill] (17.00,11.00) circle [radius=0.05];
\draw [fill] (16.86,11.14) circle [radius=0.05];
\draw [fill] ( 1.00,12.00) circle [radius=0.05];
\draw [fill] ( 5.07,11.93) circle [radius=0.05];
\draw [fill] ( 4.93,12.07) circle [radius=0.05];
\draw [fill] ( 9.14,11.86) circle [radius=0.05];
\draw [fill] ( 9.00,12.00) circle [radius=0.05];
\draw [fill] ( 8.86,12.14) circle [radius=0.05];
\draw [fill] (13.14,11.86) circle [radius=0.05];
\draw [fill] (13.00,12.00) circle [radius=0.05];
\draw [fill] (12.86,12.14) circle [radius=0.05];
\draw [fill] (17.14,11.86) circle [radius=0.05];
\draw [fill] (17.00,12.00) circle [radius=0.05];
\draw [fill] (16.86,12.14) circle [radius=0.05];

\draw ( 1.00, 1.00) --( 1.00, 0.00);
\draw ( 2.00, 1.00) --( 1.00, 0.00);
\draw ( 3.00, 1.00) --( 2.00, 0.00);
\draw ( 9.00, 1.00) --( 9.00, 0.00);
\draw (10.00, 1.00) --( 9.00, 0.00);
\draw ( 1.00, 2.00) --( 1.00, 1.00);
\draw ( 3.00, 2.00) --( 2.00, 1.00);
\draw ( 9.00, 2.00) --( 9.00, 1.00);
\draw (11.00, 2.00) --(10.00, 1.00);
\draw ( 1.00, 3.00) --( 1.00, 2.00);
\draw ( 4.93, 3.07) --( 5.00, 2.00);
\draw ( 8.93, 3.07) --( 9.00, 2.00);
\draw ( 1.00, 4.00) --( 1.00, 3.00);
\draw ( 5.07, 3.93) --( 5.07, 2.93);
\draw ( 4.93, 4.07) --( 4.93, 3.07);
\draw ( 9.00, 4.00) --( 9.07, 2.93);
\draw ( 8.86, 4.14) --( 8.93, 3.07);
\draw (10.00, 4.00) --( 9.07, 2.93);
\draw (13.00, 4.00) --(13.00, 3.00);
\draw ( 1.00, 5.00) --( 1.00, 4.00);
\draw ( 5.07, 4.93) --( 5.07, 3.93);
\draw ( 4.93, 5.07) --( 4.93, 4.07);
\draw ( 9.14, 4.86) --( 9.14, 3.86);
\draw ( 9.00, 5.00) --( 9.00, 4.00);
\draw ( 8.86, 5.14) --( 8.86, 4.14);
\draw (10.00, 5.00) --( 9.14, 3.86);
\draw (11.00, 5.00) --(10.00, 4.00);
\draw (13.00, 5.00) --(13.00, 4.00);
\draw (17.00, 5.00) --(17.00, 4.00);
\draw (18.00, 5.00) --(17.00, 4.00);
\draw ( 1.00, 6.00) --( 1.00, 5.00);
\draw ( 5.07, 5.93) --( 5.07, 4.93);
\draw ( 4.93, 6.07) --( 4.93, 5.07);
\draw ( 9.14, 5.86) --( 9.14, 4.86);
\draw ( 9.00, 6.00) --( 9.00, 5.00);
\draw ( 8.86, 6.14) --( 8.86, 5.14);
\draw (11.00, 6.00) --(10.00, 5.00);
\draw (12.93, 6.07) --(13.00, 5.00);
\draw (17.00, 6.00) --(17.00, 5.00);
\draw (19.00, 6.00) --(18.00, 5.00);
\draw ( 1.00, 7.00) --( 1.00, 6.00);
\draw ( 5.07, 6.93) --( 5.07, 5.93);
\draw ( 4.93, 7.07) --( 4.93, 6.07);
\draw ( 9.14, 6.86) --( 9.14, 5.86);
\draw ( 9.00, 7.00) --( 9.00, 6.00);
\draw ( 8.86, 7.14) --( 8.86, 6.14);
\draw (13.00, 7.00) --(13.07, 5.93);
\draw (12.86, 7.14) --(12.93, 6.07);
\draw (16.93, 7.07) --(17.00, 6.00);
\draw ( 1.00, 8.00) --( 1.00, 7.00);
\draw ( 5.07, 7.93) --( 5.07, 6.93);
\draw ( 4.93, 8.07) --( 4.93, 7.07);
\draw ( 9.14, 7.86) --( 9.14, 6.86);
\draw ( 9.00, 8.00) --( 9.00, 7.00);
\draw ( 8.86, 8.14) --( 8.86, 7.14);
\draw (13.14, 7.86) --(13.14, 6.86);
\draw (13.00, 8.00) --(13.00, 7.00);
\draw (12.86, 8.14) --(12.86, 7.14);
\draw (17.00, 8.00) --(17.07, 6.93);
\draw (16.86, 8.14) --(16.93, 7.07);
\draw (18.00, 8.00) --(17.07, 6.93);
\draw ( 1.00, 9.00) --( 1.00, 8.00);
\draw ( 5.07, 8.93) --( 5.07, 7.93);
\draw ( 4.93, 9.07) --( 4.93, 8.07);
\draw ( 9.14, 8.86) --( 9.14, 7.86);
\draw ( 9.00, 9.00) --( 9.00, 8.00);
\draw ( 8.86, 9.14) --( 8.86, 8.14);
\draw (13.14, 8.86) --(13.14, 7.86);
\draw (13.00, 9.00) --(13.00, 8.00);
\draw (12.86, 9.14) --(12.86, 8.14);
\draw (17.14, 8.86) --(17.14, 7.86);
\draw (17.00, 9.00) --(17.00, 8.00);
\draw (16.86, 9.14) --(16.86, 8.14);
\draw (18.00, 9.00) --(17.14, 7.86);
\draw (19.00, 9.00) --(18.00, 8.00);
\draw ( 1.00,10.00) --( 1.00, 9.00);
\draw ( 5.07, 9.93) --( 5.07, 8.93);
\draw ( 4.93,10.07) --( 4.93, 9.07);
\draw ( 9.14, 9.86) --( 9.14, 8.86);
\draw ( 9.00,10.00) --( 9.00, 9.00);
\draw ( 8.86,10.14) --( 8.86, 9.14);
\draw (13.14, 9.86) --(13.14, 8.86);
\draw (13.00,10.00) --(13.00, 9.00);
\draw (12.86,10.14) --(12.86, 9.14);
\draw (17.14, 9.86) --(17.14, 8.86);
\draw (17.00,10.00) --(17.00, 9.00);
\draw (16.86,10.14) --(16.86, 9.14);
\draw (19.00,10.00) --(18.00, 9.00);
\draw ( 1.00,11.00) --( 1.00,10.00);
\draw ( 5.07,10.93) --( 5.07, 9.93);
\draw ( 4.93,11.07) --( 4.93,10.07);
\draw ( 9.14,10.86) --( 9.14, 9.86);
\draw ( 9.00,11.00) --( 9.00,10.00);
\draw ( 8.86,11.14) --( 8.86,10.14);
\draw (13.14,10.86) --(13.14, 9.86);
\draw (13.00,11.00) --(13.00,10.00);
\draw (12.86,11.14) --(12.86,10.14);
\draw (17.14,10.86) --(17.14, 9.86);
\draw (17.00,11.00) --(17.00,10.00);
\draw (16.86,11.14) --(16.86,10.14);
\draw ( 1.00,12.00) --( 1.00,11.00);
\draw ( 5.07,11.93) --( 5.07,10.93);
\draw ( 4.93,12.07) --( 4.93,11.07);
\draw ( 9.14,11.86) --( 9.14,10.86);
\draw ( 9.00,12.00) --( 9.00,11.00);
\draw ( 8.86,12.14) --( 8.86,11.14);
\draw (13.14,11.86) --(13.14,10.86);
\draw (13.00,12.00) --(13.00,11.00);
\draw (12.86,12.14) --(12.86,11.14);
\draw (17.14,11.86) --(17.14,10.86);
\draw (17.00,12.00) --(17.00,11.00);
\draw (16.86,12.14) --(16.86,11.14);
\draw ( 1.00,13.00) --( 1.00,12.00);
\draw ( 5.07,12.93) --( 5.07,11.93);
\draw ( 4.93,13.07) --( 4.93,12.07);
\draw ( 9.14,12.86) --( 9.14,11.86);
\draw ( 9.00,13.00) --( 9.00,12.00);
\draw ( 8.86,13.14) --( 8.86,12.14);
\draw (13.14,12.86) --(13.14,11.86);
\draw (13.00,13.00) --(13.00,12.00);
\draw (12.86,13.14) --(12.86,12.14);
\draw (17.14,12.86) --(17.14,11.86);
\draw (17.00,13.00) --(17.00,12.00);
\draw (16.86,13.14) --(16.86,12.14);

\end{tikzpicture}
\caption{$\Ext_{\SA(1)}^{s,t}(H^*L_0, \FF_2)$} \label{fig:ExtA1L0}
\end{figure}
%
%
%
\bibliographystyle{gtart}
\bibliography{bailey_ricka}


\end{document}